\documentclass[12pt]{amsart}
\usepackage{amssymb,amsmath, amsthm, comment, multirow, tikz, cite, booktabs,
adjustbox, graphicx, subcaption, multicol, placeins, url}

\newcommand{\Z}{\mathbb{Z}}

\newcommand{\C}{\mathbb{C}}
\newcommand{\Q}{\mathbb{Q}}
\newcommand{\R}{\mathbb{R}}
\newcommand{\N}{\mathbb{N}}
\newcommand{\HH}{\mathbb{H}}

\newtheorem{theorem}{Theorem}[section]
\newtheorem{theorem*}{Theorem}
\newtheorem{algorithm}[theorem]{Algorithm}

\numberwithin{equation}{section}

\title[Computing bases of modular forms]{Computing bases of modular forms using the graded algebra structure}
\author[M.~Lam]{Michael O.~Lam}
\address{James Madison University, Department of Computer Science, MSC 4103, Harrisonburg, VA, 22807}
\email{lam2mo@jmu.edu}

\author[N.~McClelland]{Noah S.~McClelland} 

\author[M.~Petty]{Matthew R.~Petty}

\author[J.~Webb]{John J.~B.~Webb}
\address{James Madison University, Department of Mathematics and Statistics, MSC 1911, Harrisonburg, VA, 22807}
\email{webbjj@jmu.edu}

\begin{document}

\date{\today}


\date{\today}

\begin{abstract}
We develop a new algorithm to compute a basis for $M_k(\Gamma_0(N))$, the space of weight $k$ holomorphic modular forms on $\Gamma_0(N)$, in the case when the graded algebra of modular forms over $\Gamma_0(N)$ is generated at weight two.  Our tests show that this algorithm significantly outperforms a commonly used algorithm which relies more heavily on modular symbols. 
\end{abstract}

\maketitle

\section{Introduction}
\label{sec:intro}

Let $N$ be a natural number and let $k$ be an even natural number. A weight $k$ level $N$ modular form is a function $f(z)$ which is holomorphic on the complex upper-half plane $\HH = \{z=x+iy: \, x,y \in \R, \, y>0\}$ and as it approaches the cusps, $\Q \cup \{i\infty\}$, that satisfies the transformation law 
\begin{equation}\label{TransformationLaw} f\left(\dfrac{az+b}{cz+d}\right)= (cz+d)^k f(z) \end{equation}
for all matrices $\begin{pmatrix} a & b \\ c & d \end{pmatrix} \in \Gamma_0(N)$ where
\[ \Gamma_0(N) = \left\{\begin{pmatrix} a & b \\ c & d\end{pmatrix} : \, a,b,c,d\in \Z, \, N\mid c, \, ad-bc=1\right\}\, .\]
Modular forms play an important role in modern number theory---they have been used to study a wide array of objects such as quadratic forms, partitions, elliptic curves, and $L$-functions like the Riemann zeta-function.  Because $\begin{pmatrix} 1& 1\\0&1\end{pmatrix} \in \Gamma_0(N)$, \eqref{TransformationLaw} tells us that each weight $k$ level $N$ modular form is 1-periodic.  If we set $q=e^{2\pi i z}$, we can write a modular form as a Fourier series,
\[ f(z) = \sum_{n\geq 0} a_f(n) q^n \, ,\]
where each Fourier coefficient $a_f(n) \in \C$; often for a modular form these Fourier coefficients are all rational numbers, or even integers, and contain valuable combinatorial data.

We denote the set of all weight $k$ level $N$ modular forms by $M_k(\Gamma_0(N))$---this set is a finite dimensional vector space over $\C$.  Define
\begin{equation}\label{coeff_bound} B(N,k) = \frac{Nk}{12} \prod_{p \text{ prime, } p\mid N} \left(1+\frac1{p}\right) \, ,\end{equation}
then $\dim(M_k(\Gamma_0(N))) \leq B(N,k)+1$, and each modular form $f(z)\in M_k(\Gamma_0(N))$ can be uniquely identified by its Fourier coefficients $a_f(n)$ for $n\leq B(N,k)$.  (For an exact formula for $\dim(M_k(\Gamma_0(N)))$ see, for example, Theorem 3.5.1 in \cite{DS}.)  If $f(z)\in M_k(\Gamma_0(N))$ and $g(z)\in M_\ell (\Gamma_0(N))$ are two modular forms with the same level and possibly different weights, it follows from \eqref{TransformationLaw} that $f(z)g(z) \in M_{k+\ell} (\Gamma_0(N))$ and $f(z) \in M_k (\Gamma_0(tN))$ for any positive integer $t$.

Because interesting results are often obtained by studying what happens to modular forms when certain linear transformations are applied to them, it becomes important to have an accessible basis for the relevant space.  The LMFDB (\url{www.lmfdb.org}) is online database of automorphic objects and describes bases of modular forms for a limited number of weights and levels.  While there are algorithms available to find a basis for an arbitrary space of modular forms, the calculation is computationally expensive, particularly if the weight of the space is large.

This paper presents a new algorithm for calculating bases of weight $k$ modular forms on $\Gamma_0(N)$ where $N$ is composite and $\Gamma_0(N)$ has no elliptic points.  We say $z \in \HH$ is an elliptic point of $\Gamma_0(N)$ if there exists some $\begin{pmatrix} a & b \\ c & d\end{pmatrix} \in \Gamma_0(N)$, with $\begin{pmatrix} a & b \\ c & d\end{pmatrix}\neq \pm \begin{pmatrix} 1 & 0 \\ 0 & 1 \end{pmatrix}$, such that $\dfrac{az +b}{cz+d} = z$.  Whereas the most common algorithms for calculating such bases heavily relies on related objects known as modular symbols, our technique utilizes the structure of the graded ring of level $N$ modular forms,
\[ \mathcal{M}(N) = \bigcup_{k \in 2\N} M_k(\Gamma_0(N)) \, .\]
While it has long been known that $\mathcal{M}(1)$ is generated by level 1 modular forms of weight 4 and 6,  recent work of Rustom \cite{Rustom} and Voight and Zureick-Brown \cite{VDZB} shows that $\mathcal{M}(N)$ is generated by level $N$ modular forms of weight $k \leq 6$ for all $N$.  In the case where $N$ is composite and $\Gamma_0(N)$ has no elliptic points, then $\mathcal{M}(N)$ is generated entirely by weight two modular forms.  This happens precisely for composite $N$, which will be referred to as ``good'' $N$ throughout this paper, that satisfy both of the following conditions:
\begin{enumerate}
\item Either $4\mid N$ or $p \mid N$ for some prime $p \equiv 3 \pmod{4}$, 
\item Either $9\mid N$ or $p \mid N$ for some prime $p \equiv 2 \pmod{3}$.
\end{enumerate}
We note that for any $N$, at least one of $N$, $2N$, $3N$, or $4N$ must satisfy the above conditions.  This implies the efficacy of the following.

\begin{algorithm}\label{MainAlgorithm} Given $k\in 2\N$ and a good $N$, a basis for the space $M_k(\Gamma_0(N))$ is calculated by
\begin{enumerate}
\item Calculating a basis $\mathcal{B}_2$ for $M_2(\Gamma_0(N))$ and
\item Calculating products of $k/2$ members of $\mathcal{B}_2$ until $\dim(M_k(\Gamma_0(N))$ linearly independent modular forms are found.
\end{enumerate}
\end{algorithm}

In the next section, we compare the performance of this algorithm against a standard algorithm that relies more heavily on modular symbols.  In particular, our computations show significant improvement in both speed and memory usage, particularly for spaces of high weight.  In Section 3, we provide a detailed description of our implementation of Algorithm \ref{MainAlgorithm}, while in Section 4, we investigate building bases of modular forms using a specific kind of modular form known as eta-quotients.

\section{Performance of Algorithm \ref{MainAlgorithm}}

We implemented Algorithm \ref{MainAlgorithm} using the open source computer algebra system SageMath (referred to here simply as Sage).  This system provides extensive support for modular forms, including the ability to calculate bases of modular forms for a wide array of congruence subgroups.  We compared the performance of our implementation of Algorithm \ref{MainAlgorithm} against the built-in Sage algorithm for good $N \leq 198$ and a range of weights.  Sage's internal algorithm uses modular symbols to calculate a substantial portion of each basis and is very similar to the algorithm used by Magma, another computer algebra system.  Further details on Sage's internal algorithm are in \cite{Stein}.  As $N$ and $k$ increase, both the time and memory required to complete the basis computation increase significantly, which limited the number of weights for which we were able to test the two algorithms for each good $N$.  We will present here a portion of the data the tests produced, which are typical of the overall results.  The full data set can be found at \url{http://educ.jmu.edu/~webbjj/modformcalc}.  Our experiments were run
on a 24-core AMD Opteron 6344 server with 64GB of RAM running Sage 6.4.1 in CentOS 6.8. 

\begin{table} 
\centering
\begin{tabular}{ c | r | r | r | r }
Weight & \parbox[c]{2cm}{\center{Sage\\Time (sec)}} & \parbox[c]{2cm}{\center{Alg 1.1\\Time (sec)}} & \parbox[c]{2cm}{\center{Sage\\Memory (MB)}} & \parbox[c]{2cm}{\center{Alg 1.1\\Memory (MB)}} \\
\hline
12 & 0.7 & 0.2 & 16.6 & 8.8 \\
24 & 5.2 & 0.3 & 37.4 & 8.9 \\
36 & 32.9 & 0.6 & 76.3 & 9.0 \\
48 & 161.7 & 1.0 & 140.0 & 9.2 \\
60 & 732.2 & 2.0 & 332.8 & 9.4 \\
72 & 2171.7 & 2.8 & 672.8 & 9.6 \\
84 & 6656.9 & 4.5 & 1394.9 & 9.9 \\
96 & 15073.1 & 7.2 & 2640.8 & 10.3 \\
\hline  

\end{tabular}  

\caption{Level 8 Time and Memory Allocation Comparisons}\label{lev8}   
  
\end{table}

\subsection*{Level 8}
For level 8, both algorithms were tested for all even weights $k \leq 100$.  In Table \ref{lev8}, we focus on weights that are multiples of 12.  The time for Sage's algorithm to complete roughly tripled and the memory allocation doubled for each increase in the weight by 12, while Algorithm \ref{MainAlgorithm}'s time at most doubled, and the memory allocation increased linearly.  At weight 96, which has dimension 97, Sage's algorithm took over 4 hours to compute a basis, while Algorithm \ref{MainAlgorithm} computed the basis in just over 7 seconds, with a smaller memory allocation by a factor of over 250.

\begin{table}
\centering
\begin{tabular}{ c | r | r | r | r }
Weight & \parbox[c]{2cm}{\center{Sage\\Time (sec)}} & \parbox[c]{2cm}{\center{Alg 1.1\\Time (sec)}} & \parbox[c]{2cm}{\center{Sage\\Memory (MB)}} & \parbox[c]{2cm}{\center{Alg 1.1\\Memory (MB)}} \\
\hline
4 & 1.2 & 0.6 & 15.5 & 11.6 \\
8 & 14.2 & 1.0 & 97.0 & 12.4 \\
12 & 98.6 & 2.0 & 227.5 & 13.1 \\
16 & 566.9 & 3.0 & 501.8 & 13.8 \\
20 & 2071.7 & 10.5 & 1326.3 & 14.5 \\
24 & 12727.2 & 23.3 & 2972.3 & 15.3 \\
\hline  

\end{tabular}  

\caption{Level 36 Time and Memory Allocation Comparisons}\label{lev36}   
  
\end{table}

\subsection*{Level 36}
For level 36, Sage's algorithm was tested for all even weights up to 24, while Algorithm \ref{MainAlgorithm} was tested up to weight 62.  In Table \ref{lev36}, we focus on weights that are multiples of 4.  At weight 24, which has dimension 144, Algorithm \ref{MainAlgorithm} calculates the basis faster by a factor of over 500, while using substantially less memory.  At weight 62, which has dimension 372, it does take Algorithm \ref{MainAlgorithm} over 6 hours to calculate a basis, but only uses 23.7 MB of memory. 
 
\begin{table}
\centering
\begin{tabular}{ c | r | r | r | r }
Weight & \parbox[c]{2cm}{\center{Sage\\Time (sec)}} & \parbox[c]{2cm}{\center{Alg 1.1\\Time (sec)}} & \parbox[c]{2cm}{\center{Sage\\Memory (MB)}} & \parbox[c]{2cm}{\center{Alg 1.1\\Memory (MB)}} \\
\hline
4 & 13.7 & 2.1 & 65.0 & 21.0 \\
6 & 76.4 & 4.5 & 527.2 & 25.1 \\
8 & 324.9 & 12 & 730.1 & 29.8 \\
10 & 1097.2 & 32.3 & 1963.3 & 38.3 \\
\hline  

\end{tabular}  

\caption{Level 105 Time and Memory Allocation Comparisons} \label{lev105}  
  
\end{table}

\subsection*{Level 105}
At level 105, Sage's algorithm was tested up to weight 10, while Algorithm \ref{MainAlgorithm} computed bases up to weight 30.  We omit the data for weight 2 when both algorithms are essentially identical.  In Table \ref{lev105}, we see again that both the time and memory allocation increase at a much faster rate than those for Algorithm \ref{MainAlgorithm}.  At weight 30 where the dimension is 468, Algorithm \ref{MainAlgorithm} computed a basis in 36 hours and uses only 364.2 MB of memory.

\begin{table}
\centering
\begin{tabular}{ c | r | r | r | r }
Weight & \parbox[c]{2cm}{\center{Sage\\Time (sec)}} & \parbox[c]{2cm}{\center{Alg 1.1\\Time (sec)}} & \parbox[c]{2cm}{\center{Sage\\Memory (MB)}} & \parbox[c]{2cm}{\center{Alg 1.1\\Memory (MB)}} \\
\hline
4 & 245 & 14.9 & 531.7 & 79.0 \\
6 & 2196.4 & 29.4 & 2915.0 & 107.4 \\
8 & 12922.0 & 54.1 & 10398.4 & 139.3 \\
10 & 30895.8 & 102.3 & 25693.9 & 170.0 \\
\hline  

\end{tabular}  

\caption{Level 198 Time and Memory Allocation Comparisons} \label{lev198}  
  
\end{table}

\subsection*{Level 198}
The highest level we compared the two algorithms was 198, where Sage's algorithm was tested up to weight 10, and Algorithm \ref{MainAlgorithm} was tested up to weight 18.  Table \ref{lev198} shows Algorithm \ref{MainAlgorithm} performed significantly better again.  At weight 18 where the dimension is 620, Algorithm \ref{MainAlgorithm} computed a basis in 2.5 hours, using 350.7 MB of memory.

\section{Implementation of Algorithm \ref{MainAlgorithm}}
Our source code for the implementation of Algorithm \ref{MainAlgorithm} is available at \url{http://educ.jmu.edu/~webbjj/modformcalc/main.sage}.  While Sage has robust support for modular forms, we found it necessary to define a new object class for modular forms which involves the graded algebra structure.  In particular, for a given basis of $M_2(\Gamma_0(N))$, we save the representation of each modular form as a linear combination of products from this basis as well as the Fourier expansion.  This representation is useful if it is necessary to calculate the Fourier coefficients of the forms to a much higher degree, such as when one is calculating the action of Hecke operators on the space.  In this case, it suffices to calculate the additional Fourier coefficients of the weight two basis forms and then use this representation to quickly lift the accuracy to a higher degree.  All of our calculations are performed over $\Q$, although any other fields could be substituted without needing to make significant changes.  

Algorithm \ref{MainAlgorithm} begins by calculating a basis for $M_2(\Gamma_0(N))$.  The authors found that Sage's internal algorithm for computing bases of modular forms is extremely efficient at weight two, although in Section 4 an alternate method which computes a basis of weight two eta-quotients is detailed.  With this basis in hand, we now start calculating forms in our desired space, $M_k(\Gamma_0(N))$.  At all times, our saved set of linearly independent weight $k$ forms is upper triangular, which we now describe.  For a modular form $f(z) = \sum_{n\geq 0} a_f(n) q^n$, we let $v_\infty (f)= \min \{n : a_f(n) \neq 0\}$ denote the degree of vanishing (also called the degree of the zero) of $f$ at the cusp $i\infty$.  We say that a set $\{g_1, \ldots, g_t\}$ is upper triangular if $v_\infty(g_1)<v_\infty(g_2) < \ldots < v_\infty (g_t)$, that is, the degrees of vanishing at $i\infty$ of the basis forms are strictly increasing. As with upper triangular matrices generally in linear algebra, upper triangular bases of modular forms often are very efficient in calculations.  In particular, the following simple iterative procedure determines if another form $f$ is a linear combination of $g_1, \ldots, g_t$ and, if not, finds a form to adjoin to the set that keeps it upper triangular.

\begin{algorithm}\label{upper triangular alg} 
The upper triangular set $\mathcal{B}=\{g_1, \ldots, g_t\}$ and $f$ are taken as input.  Let $g_i = \sum_{n\geq n_i} a_{g_i}(n)q^n$ where $n_i = v_\infty (g_i)$ and $f= \sum_{n\geq n_0} a_f(n) q^n$.  For $1\leq i \leq t$, check $v_\infty (g_i)$ and $v_\infty (f)$:
\begin{itemize}
\item If $v_\infty (g_i) = v_\infty (f)$ then 
\begin{itemize}
\item Replace $f$ with $f-\tfrac{a_f(n_i)}{a_{g_i} (n_i)} g_i$.  
\item If this new $f =0$, then the original $f$ is in the span of $\mathcal{B}$ and $\mathcal{B}$ is returned unchanged.  If $f\neq 0$, increase $i$ to $i+1$.
\end{itemize}
\item Else if $v_\infty (g_i) > v_\infty (f)$ then 
\begin{itemize}
\item Insert $f$ into the $i$-th position of $\mathcal{B}$, shifting the index of $g_i, \ldots, g_t$ up by 1, and
\item Return $\mathcal{B}$ and stop the algorithm.
\end{itemize}
\item Else if $v_\infty (g_i) < v_\infty (f)$ then increase $i$ to $i+1$.
\end{itemize}
After the check for $i=t$, if $f \neq 0$, then $f$ is inserted into the $t+1$-th position of $\mathcal{B}$.  This larger $\mathcal{B}$ is returned and the algorithm ends.
\end{algorithm}

To cut down on the number of multiplications needed, we save some information from the intermediate steps.  So, for example, if we let $\{f_1, f_2, \ldots, f_t\}$ be the weight 2 basis, the first weight $k$ form computed will be $f_1^{k/2}$.  This form becomes the first member of $\mathcal{B}$, our set of linearly independent weight $k$ forms.  We save $f_1^2, f_1^3, \ldots f_1^{k/2-1}$ along the way, because then $f_1^{k/2-1} f_2, f_1^{k/2-1}f_3, \ldots, f_1^{k/2} f_t$ can then be computed with only one additional multiplication for each.  Now $f_1^{k/2-1}$ can be discarded, $f_1^{k/2-2}$ will be multiplied by $f_2$, this is saved, and then $f_1^{k/2-2}f_2^2, f_1^{k/2-2}f_2f_3, \ldots, f_1^{k/2-2}f_2f_t$ are each computed.  The search continues in this manner.  After each new weight $k$ form that is calculated, we check if it is in the span of $\mathcal{B}$, and if not, adjoin a suitable form to $\mathcal{B}$ that keeps it upper triangular.  Once $\mathcal{B}$ contains $\dim M_k(\Gamma_0(N))$ forms, it is a basis for the space and the search concludes.  

In our initial tests, we found that the order of the forms in the weight two basis could significantly affect the speed of the search.  If the initial weight two basis was upper triangular, then we have
\[ v_\infty (f_1^{k/2})< v_\infty(f_1^{k/2-1}f_2) < \ldots < v_\infty(f_1^{k/2-1}f_t) \]
so they together are linearly independent and upper triangular---the first $t$ elements of $\mathcal{B}$ have been found.  While it is quick to take the given weight two basis and produce an upper triangular basis out of it, this complicates the representation of the form in terms of the original basis which is being saved along with Fourier expansion of each form; in turn, this bogs down the calculation of the higher weight forms' representations.  Looking closely at Algorithm \ref{upper triangular alg}, if $\{g_1, \ldots, g_t\}$ is the upper triangular basis produced by iterating this algorithm on the original weight 2 basis, then there is a re-ordering of this original basis $\{f_{\sigma(1)}, f_{\sigma(2)}, \ldots , f_\sigma(t)\}$ such that $g_1 = f_{\sigma(1)}$ and $g_i$ is in the span of $f_{\sigma(1)}, \ldots , f_{\sigma(i)}$ for all $i>1$.  This implies that for each $i$, the algebra generated by $f_{\sigma(i)}, \ldots , f_{\sigma(i)}$ is equal to that generated by $g_1, \ldots , g_i$.  We use this reordered basis in our implementation because it guarantees that many forms in our basis will be found at the beginning of the search, while keeping the the representation of the forms in terms of the original weight two basis relatively simple.  Utilizing this idea of finding forms with different degrees of vanishing at $i\infty$, we simultaneously run the search starting at different places.  To be precise, while we start with $f_{\sigma(1)}^{k/2}$ and start swapping higher indexed forms into the product, we also look at $f_{\sigma(t)}^{k/2}$ and swap in lower indexed forms, and even with $f_{\sigma(\lfloor t/2 \rfloor)}^{k/2}$, where we make two more branches of the search by swapping in higher or lower indexed forms.

We note that this search algorithm we employed is straight-forward.  Potentially, this can optimized much further.  However, as seen in the data of the previous section, this implementation shows tremendous improvements over the current algorithm used to calculate bases of modular forms.
 
\section{Building bases with eta-quotients}
Dedekind's eta function, $\eta(z)$ is defined as
\begin{equation}\label{etafunc}
\eta(z) = q^{\frac1{24}}\prod_{n\geq1} (1-q^n) = q^{\frac1{24}}(1+\sum_{n\geq1} (-1)^n(q^{\frac{n(3n+1)}{2}} + q^{\frac{n(3n-1)}{2}}) \, .
\end{equation}
This function plays an important role in number theory.  For example, $1/\eta(z)$ is a generating function for partitions and $\eta(z)^{24}=\Delta(z)$ is Ramanujan's Delta function, also known as the modular discriminant function.  

For some $N>0$, we say a function $f(z)$ is a level $N$ eta-quotient if $f(z) = \prod_{\delta \mid N} \eta(\delta z)^r_\delta$, where each $r_\delta \in \Z$.  Such an $f(z)$ is a modular form in $M_k(\Gamma_0 (N))$ precisely when 
\begin{itemize}
\item it is holomorphic at the cusps, 
\item the sums $ \sum_{\delta \mid N} \delta r_\delta$, $\sum_{\delta \mid N} \frac{N}{\delta} r_\delta$ are both divisible by 24, and 
\item $\prod_{\delta \mid N} \delta^{r_\delta}$ is a rational square, 
\end{itemize}
where the weight $k = \frac12 \sum_{\delta \mid N} r_\delta$.  See \cite{Ono}, Section 1.4, for further details.   

The representation of $\eta(z)$ as a sum in \eqref{etafunc} means that the calculation of Fourier expansions of eta-quotients is computationally straightforward.  For this reason, it is natural to ask when a modular form can be represented in terms of eta-quotients.  In \cite{RouseWebb}, Rouse and the fourth author found that the entire graded algebra of level $N$ modular forms, $\mathcal{M}(N)$, is generated by level $N$ eta-quotients if and only if $M_2(\Gamma_0(N))$ is generated by eta-quotients.  While this implies that $\Gamma_0(N)$ does not have elliptic points, this is not a sufficient condition for $\mathcal{M}(N)$ to be generated by eta-quotients.  For example, $\Gamma_0(68)$ does not have elliptic points and, hence, $\mathcal{M}(N)$ is generated by weight two forms, but $M_2(\Gamma_0(68)$ is not generated by level $N$ eta-quotients. There are precisely 121 levels $N \leq 500$, identified in \cite{RouseWebb}, such that $\mathcal{M}(N)$ is generated by level $N$ eta-quotients.  To find such $N$, Rouse and the fourth author enumerated all of the eta-quotients in $M_2(\Gamma_0(N))$ by finding which tuples of exponents $(r_\delta)_{\delta\mid N}$ satisfied the conditions above, and then checking the dimension of the space spanned by these forms.

In the current work, we calculated eta-quotients using the the possible orders of vanishing these forms can have at the cusps to build bases.  In more technical terms, we found eta-quotients in $M_2(\Gamma_0(N))$ by their divisors on the modular curve $X_0(N)$.  From \eqref{etafunc}, $\eta(z)$ is non-zero on $\HH$, so all of the zeros (and poles) of an eta-quotient must occur at the cusps.  The following result of Ligozat \cite{Ligozat} calculates the order of vanishing of an eta-quotient at a rational cusp.

\begin{theorem*}
Let $c$, $d$ and $N$ be positive integers with $d | N$ and $\gcd(c,d) = 1$.
If $f(z)$ is an eta-quotient, then the order of vanishing of $f(z)$ at
the cusp $c/d$ is
\[
  \frac{N}{24} \sum_{\delta | N} \frac{\gcd(d,\delta)^{2} r_{\delta}}{\gcd(d,N/d) d \delta}\, .
\]
\end{theorem*}

We note that the order of vanishing at $c/d$ depends on the denominator $d$, but not on the numerator $c$, and the total number of zeros of a weight $k$ modular form (including multiplicities) is $B(N,k)$ from \eqref{coeff_bound}.  These relations create a linear transformation from tuple of exponents $(r_\delta)_{\delta\mid N}$ to the divisor of the eta-quotient---this transformation must be invertible because divisors of modular forms are unique up to constant multiples and all eta quotients have 1 as their first non-zero coefficient.  To find eta-quotients in $M_2(\Gamma_0(N))$, we inverted the $d(N)$ by $d(N)$ matrix generated by Ligozat's relations, where $d(N)$ is the number of divisors of $N$, and applied it to possible divisors of weight two level $N$ eta-quotients.  The advantage of this approach is that the set of possible divisors is significantly smaller than the set of possible $d(N)$ tuples of exponents (see \cite{Bhattacharya} or \cite{RouseWebb} for sharp bounds on the size of exponents for eta-quotients).

In our search for level $N$ eta-quotients to build a basis, we employ the following additional strategies.  Because $M_2(\Gamma_0(\delta)) \subset M_2(\Gamma_0(N))$ for all $\delta \mid N$, we begin our search initially in the smaller spaces.  For each eta-quotient $f(z) \in M_2(\Gamma_0(\delta))$, we use that $f(tz) \in M_2(\Gamma_0(N))$ for all $t | (N/\delta)$ will be an eta-quotient as well.  Finally, we used that if $f(z)$ is an eta-quotient in $M_2(\Gamma_0(N))$, then its image under the Fricke involution (in short, swapping the exponents $r_\delta$ and $r_{N/\delta}$ for each $\delta\mid N$) is also an eta-quotient in the space.

\begin{table} 
\centering
\begin{tabular}{ c | r | r | r }
Weight & \parbox[c]{3cm}{\center{Sage\\Time (sec)}} & \parbox[c]{3cm}{\center{Alg 1.1\\Time (sec)}} & \parbox[c]{3cm}{\center{eta-Quotient\\Time (sec)}}  \\
\hline
12 & 0.7 & 0.2 &  0.3 \\
24 & 5.2 & 0.3 &  0.3 \\
36 & 32.9 & 0.6 & 0.7 \\
48 & 161.7 & 1.0 & 1.1 \\
60 & 732.2 & 2.0 &  1.6 \\
72 & 2171.7 & 2.8 &  2.9 \\
84 & 6656.9 & 4.5 &  4.7 \\
96 & 15073.1 & 7.2 &  7.9 \\
\hline  

\end{tabular}  

\caption{Level 8 Time Comparisons with eta-Quotient Algorithm}\label{lev8eta}   
  
\end{table}

\begin{table}
\centering
\begin{tabular}{ c | r | r | r }
Weight & \parbox[c]{3cm}{\center{Sage\\Time (sec)}} & \parbox[c]{3cm}{\center{Alg 1.1\\Time (sec)}} & \parbox[c]{3cm}{\center{eta-Quotient\\Time (sec)}}  \\
\hline
4 & 1.2 & 0.6 & 1.0  \\
8 & 14.2 & 1.0 & 1.4 \\
12 & 98.6 & 2.0 & 3.7 \\
16 & 566.9 & 3.0 & 21.7 \\
20 & 2071.7 & 10.5 & 31.0 \\
24 & 12727.2 & 23.3 & 304.9 \\
\hline  

\end{tabular}

\caption{Level 36 Time Comparisons with eta-Quotient Algorithm}\label{lev36eta}   
  
\end{table}

\begin{table}
\centering
\begin{tabular}{ c | r | r | r | r }
Weight & \parbox[c]{3cm}{\center{Sage\\Time (sec)}} & \parbox[c]{3cm}{\center{Alg 1.1\\Time (sec)}} & \parbox[c]{3cm}{\center{eta-Quotient\\Time (sec)}}  \\
\hline
4 & 13.7 & 2.1 & 545.6\\
6 & 76.4 & 4.5 & 612.8 \\
8 & 324.9 & 12 & 1326.0 \\
10 & 1097.2 & 32.3 & 8815.8 \\
\hline  

\end{tabular}  

\caption{Level 105 Time Comparisons with eta-Quotient Algorithm} \label{lev105eta}  
  
\end{table}

Tables \ref{lev8eta} -- \ref{lev105eta} show some time comparisons of using eta-quotients to compute bases for higher weight spaces.  We omit the memory usage information because the amount of memory used by eta-quotient algorithm and Algorithm \ref{MainAlgorithm} is very similar.  Our data shows that for small $N$, the eta-quotient algorithm outperforms Sage's internal algorithm, but is generally slower than Algorithm \ref{MainAlgorithm}.  As $N$ increases, and the total number of eta-quotients in $M_2(\Gamma_0(N))$ grows, the performance of the eta-quotient algorithm declines significantly.    

\section{Conclusions and future work}

Our data gives clear evidence that Algorithm \ref{MainAlgorithm} outperforms the current processes that are used to compute bases of modular forms for these good levels $N$ both in terms of time and computer memory usage.  Considering that for a general level $N$, the graded algebra of modular forms $\mathcal{M}(N)$ is generated by forms of weight no more than 6, this work raises the question of whether this approach generalizes effectively to compute any space $M_k(\Gamma_0(N))$.

\section*{Acknowledgments}
The second and third authors received funding for this project through the James Madison University Tickle Scholarship Fund, while the fourth author received support from a James Madison University College of Science and Mathematics Faculty Summer Grant.

\bibliographystyle{amsplain}
\bibliography{citations}{}

\end{document}